\documentclass[a4paper,11pt]{article}
\usepackage{latexsym,amsfonts,amsmath}
\usepackage{amssymb}
\usepackage{color,soul} %,ulem}
\usepackage{xcolor}
\usepackage[mathscr]{eucal}
\usepackage{enumerate}
\usepackage{enumitem}

\usepackage{pgf}
\usepackage{tikz}
\usetikzlibrary{arrows,automata}

\usepackage{fullpage}

\definecolor{Dgreen}{rgb}{102,255,102}
\definecolor{Purp}{rgb}{102,0,204}

\def\1{\mathbf{1}}
\def\0{\mathbf{0}}

\def\NN{\mathbb{N}}

\def\RR{\mathbb{R}}

\def\PP{\Pp}

\def\car{{\mathcal{C}ar}}

\def\Pp{\boldsymbol{\mathcal{P}}}

\def\car{{\mathcal{C}ar}}

\def\XX{\mathbf{X}}
\def\AA{\mathbf{A}}
\def\KK{\mathbf{K}}
\def\SS{\mathbf{S}}
\def\HH{\mathbf{H}}
\def\MM{\boldsymbol{\mathcal{M}}}

\newcommand{\widebar}[1]{\overline{#1}}

\newcommand{\bcal}[1]{\boldsymbol{\mathcal{#1}}}

\newcommand{\bfrak}[1]{\boldsymbol{\mathfrak{#1}}}

\def\ds{\displaystyle}

\definecolor{Dgreen}{rgb}{0,0.5,0}

\usepackage{natbib}
\newcounter{hypot}

\newcounter{assump}
\renewcommand{\theassump}{$(\mathrm{\Alph{hypot}}_\arabic{assump})$}

\newcounter{condS}
\renewcommand{\thecondS}{$(\mathrm{S}_\arabic{condS})$}
\newenvironment{condit}{\begin{list}
      {\hspace{\labelsep} \bfseries \thecondS}
      {\leftmargin=0pt
       \labelwidth=0cm
       \usecounter{condS}
       \def\makelabel##1{##1}
       }}{\end{list}}

\begin{document}
\newtheorem{theorem}{Theorem}[section]
\newtheorem{proposition}[theorem]{Proposition}
\newtheorem{lemma}[theorem]{Lemma}
\newtheorem{corollary}[theorem]{Corollary}
\newtheorem{definition}[theorem]{Definition}
\newtheorem{remark}[theorem]{Remark}
\newtheorem{example}[theorem]{Example}
\newtheorem{conjecture}[theorem]{Conjecture}
\newtheorem{assumption}[theorem]{Assumption}

\bibliographystyle{plain}

\title{Absorbing Markov Decision Processes
\footnote{This is a revised version of the submission {\tt https://arxiv.org/abs/2309.07059v1}\  We have weakened the Condition \ref{Hyp-Q-strongly-continuous} at the beginning of Section \ref{Sec-4} with all the results in the paper remaining valid. 
 Supported by grant PID2021-122442NB-I00 from the Spanish \textit{Ministerio de Ciencia e Innovaci\'on.}}}

\author{Fran\c{c}ois Dufour\footnote{
Institut Polytechnique de Bordeaux; INRIA Bordeaux Sud Ouest, Team: ASTRAL; IMB, Institut de Math\'ematiques de Bordeaux, Universit\'e de Bordeaux, France \tt{francois.dufour@math.u-bordeaux.fr}}
\and
Tom\'as Prieto-Rumeau\footnote{Statistics Department, UNED, Madrid, Spain. e-mail: {\tt{tprieto@ccia.uned.es}}\quad {(Author for correspondence)}}}
\maketitle
\begin{abstract}
In this paper, we study discrete-time absorbing Markov Decision Processes (MDP) with measurable state space and Borel action space with a given initial distribution.
For such models, solutions to the characteristic equation that are not occupation measures may exist.
Several necessary and sufficient conditions are provided to guarantee that any solution to the characteristic equation is an occupation measure.
Under the so-called continuity-compactness conditions, we first show that a measure is precisely an occupation measure if and only if it satisfies the characteristic equation and an additional absolute continuity condition.
Secondly, it is shown that the set of occupation measures is compact in the weak-strong topology if and only if the model is uniformly absorbing. 
Several examples are provided to illustrate our results.
\end{abstract}
{\small 
}

\section{Introduction}

In this work, we consider a discrete-time absorbing Markov Decision Process (MDP) with measurable state space $\XX$, Borel action space $\AA$, with a given initial distribution denoted by $\eta$.
We consider a general measurable state space to cover the models studied in game theory. Indeed, most models studied in nonzero-sum games theory assume that the state space is a general measurable space; 
see for example \cite{he-sun17,nowak92} and the references therein. Recently, nonzero-sum Markov games have been studied for absorbing models in \cite{dufour23}, which has  motivated the work presented in this paper.

An absorbing MDP is a special type of MDP where some measurable subset $\Delta$ of the state space is considered as \textit{absorbing} and leads to the stopping of the process.
This means that once the system enters in $\Delta$, it stays there indefinitely with zero reward or cost. The underlying assumption is that for any policy, the expected entrance time to  $\Delta$ is finite.
For more details on this type of MDP, we refer to the following papers and books \cite{altman99,piunovskiy19,feinberg12,piunovskiy13} and their references.
It can be easily shown that a discounted model can be reformulated as an absorbing model by adding an additional state and modifying the transition and performance function to make the singleton associated with that state the absorbing set. For more details on this construction, we refer to \cite[page 137]{altman99} as well as \cite[page 132]{feinberg12}. 
However, the notion of an absorbing model is more general than that of a discounted model, in the sense that there exist absorbing models that cannot be transformed into discounted models.

\bigskip

The objective of this paper is to study the properties of the solutions of the characteristic equation and of the set of occupation measures.
The notion of occupation measure is particularly important and plays a central role in solving constrained Markov decision processes (see for example the references
\cite{altman99,borkar02,hernandez96,hernandez99,piunovskiy97}). It describes the expected amount of time spent by the state and action processes
in any measurable subset of $\XX\times\AA$ under any policy. For an absorbing model, any occupation measure satisfies the so-called characteristic equation which is of the form
\begin{align}
\label{eq-charac-intro}
\mu^\XX (\cdot)= \eta_{\Delta^c}(\cdot) + \int_{\XX\times \AA} Q\mathbb{I}_{\Delta^c}(\cdot |x,a) \mu(dx,da)
\end{align}
(see equation \eqref{eq-linear-equation} for a more precise statement) where $\eta_{\Delta^c}$ is the restriction of the measure $\eta$ to~$\Delta^c$, $\mu^\XX$ is the marginal  measure of  $\mu$ on $\XX$, and $Q\mathbb{I}_{\Delta^c}$ is the product of the transition kernel $Q$ of the MDP and $\mathbb{I}_{\Delta^c}$, the indicator kernel of the set $\Delta^c$.
However, there may exist solutions to the characteristic equation \eqref{eq-charac-intro} that are not occupation measures (see example \ref{example-phantom}).
These are called phantom solutions to the characteristic equation (or simply phantom measures).
This is a notable difference from the discounted case, for which there are no phantom measures for the associated characteristic equation, which implies that a solution of the characteristic equation is necessarily an occupation measure. In order to avoid this pathological phenomenon for absorbing models, we will provide in Theorem \ref{prop-characterizations} several necessary and sufficient conditions to ensure that there are no phantom measures, thereby guaranteeing that any solution to the characteristic equation is an occupation measure.

By strengthening our hypotheses, that is, assuming the so-called continuity-compactness conditions introduced by Sch\"{a}l in \cite{schal75} (see Condition (S) below), we will show 
in our second main result (see Theorem \ref{Eq-Characteristic++}) that a measure is an occupation measure if and only if it satisfies the characteristic equation \eqref{eq-charac-intro} and an additional condition.
Let us emphasize that a measure satisfies this last condition if its marginal on $\XX$ is absolutely continuous with respect to a reference probability measure.
The construction of this reference measure is essentially based on Condition (S). An interesting property that plays a central role in the proof of our last result (see Theorem \ref{Relative-compactness-set-D}) is that if the initial distribution $\eta$ of the MDP is replaced by this reference probability measure, then the MDP remains absorbing.

Finally, it will be shown in our last main result (see Theorem \ref{Relative-compactness-set-D}) that under the continuity-compactness conditions 
the set of occupation measures is compact in the $ws$-topology (see \cite{balder01} for more details on the weak-strong topology) if and only if the model is uniformly absorbing(the notion of a uniformly absorbing MDP was introduced in \cite[Definition 3.6]{piunovskiy19}).
In fact, to be more precise, our result also shows that the set of occupation measures is relatively compact in the $ws$-topology  if and only if the model is uniformly absorbing.
In particular, it highlights the fact that the assumption that the MDP is absorbing is not sufficient to guarantee the (relative) compactness of the set of occupation measures. 
This last point is illustrated by an example. Observe that a discounted MDP is  uniformly absorbing and so we recover the classical result that, under the continuity-compactness conditions, the set of occupation measures of a discounted model is compact. The question of the compactness of the set of occupation measures is central to solving constrained optimal control problems using linear programming.

The main results of this paper appear to be more general than those established in the literature, in particular the characterization of phantom and occupation measures (Theorems \ref{prop-characterizations} and \ref{Relative-compactness-set-D}), as well as the compactness of occupation measures (Theorem \ref{Relative-compactness-set-D}), notably because these characterizations are made by means of necessary and sufficient conditions.

The rest of the paper is organized as follows. In the remaining of this section we introduce some notation and recall a standard result that will be useful in the sequel.
In Section \ref{Sec-2}, we introduce the control model under consideration and provide some basic definitions.
Section \ref{Sec-3} is devoted to the analysis of the characteristic equation.
In the last section, it is shown that the occupation measures are characterized by the characteristic equation and an additional condition. We also study the compactness property of the set of occupation measures.

\paragraph{Notation and terminology.}
The power set of a set $S$  (the set of all subsets of $S$) is denoted by~$2^{S}$.
On a measurable space $(\mathbf{\Omega},\mathcal{F})$ we will consider the set of finite signed measures $\bcal{M}(\mathbf{\Omega})$, the set of finite nonnegative measures
$\bcal{M}^+(\mathbf{\Omega})$, and the set of probability measures~$\bcal{P}(\mathbf{\Omega})$.
For a set $\Gamma\in\mathcal{F}$, we denote by $\mathbf{I}_{\Gamma}:\Omega\rightarrow\{0,1\}$ the indicator function of the set~$\Gamma$, that is,
$\mathbf{I}_{\Gamma}(\omega)=1$ if and only if $\omega\in\Gamma$.
For $\omega\in\mathbf{\Omega}$, we write $\delta_{\{\omega\}}$ for the Dirac probability measure at $\omega$ defined on $(\mathbf{\Omega},\mathcal{F})$ by
$\delta_{\{\omega\}}(B)=\mathbf{I}_{B}(\omega)$ for any $B\in\mathcal{F}$. If $\mu\in\bcal{M}(\mathbf{\Omega})$ and $\Gamma\in\mathcal{F}$, we denote by
$\mu_{\Gamma}$ the restriction of the measure $\mu$ to $\Gamma$, that is, the measure on $(\mathbf{\Omega},\mathcal{F})$ defined by $\mu_{\Gamma}(B)=\mu(\Gamma\cap B)$ for $B\in\mathcal{F}$.
The trace $\sigma$-algebra of a set $\Gamma\subseteq\mathbf{\Omega}$ is denoted by~$\mathcal{F}_{\Gamma}$.
On $\bcal{P}(\mathbf{\Omega})$, the $s$-topology is the coarsest topology that makes  $\mu\mapsto \mu(D)$ continuous for every $D\in\mathcal{F}$.

Given a measurable space  $(\mathbf{\Omega},\mathcal{F})$ and $\lambda\in\bcal{P}(\mathbf{\Omega})$, we will 
denote by $L^{1}(\mathbf{\Omega},\mathcal{F},\lambda)$  the family of measurable functions  (identifying those which are $\lambda$-a.s. equal)  $f:\mathbf{\Omega}\rightarrow\RR$ which are $\lambda$-integrable, i.e.,
$\|f\|_1=\int_\mathbf{\Omega} |f(\omega)|\lambda(d\omega)<\infty$. We note that, throughout this paper, the real numbers set $\RR$ will be always endowed with its Borel $\sigma$-algebra.

Let $(\mathbf{\Omega},\mathcal{F})$ and $(\widetilde{\mathbf{\Omega}},\widetilde{\mathcal{F}})$ be two measurable spaces.
A kernel on $\widetilde{\mathbf{\Omega}}$ given $\mathbf{\Omega}$ is a mapping
$Q:\mathbf{\Omega}\times\widetilde{\mathcal{F}}\rightarrow\RR^+$ such that $\omega\mapsto Q(B|\omega)$ is 
measurable on $(\mathbf{\Omega},\mathcal{F})$ for every $B\in\widetilde{\mathcal{F}}$,   and  $B\mapsto Q(B|\omega)$ is in $\bcal{M}^+(\widetilde{\mathbf{\Omega}})$
for every $\omega\in\mathbf{\Omega}$. If $Q(\widetilde{\mathbf{\Omega}}|\omega)=1$ for all $\omega\in\mathbf{\Omega}$ then we say that $Q$ is a \textit{stochastic kernel}.
We write $\mathbb{I}_{\Gamma}$ for the kernel on $\mathbf{\Omega}$ given $\mathbf{\Omega}$ defined by
$\mathbb{I}_{\Gamma}(B|\omega)=\mathbf{I}_{\Gamma}(\omega) \delta_{\{\omega\}}(B)$ for $\omega\in\mathbf{\Omega}$ and $B\in\mathcal{F}$.
Let $Q$ be a stochastic kernel on $\widetilde{\mathbf{\Omega}}$ given $\mathbf{\Omega}$.
For  a bounded measurable function $f:\widetilde{\mathbf{\Omega}}\rightarrow\RR$, we will denote by $Qf:\mathbf{\Omega}\rightarrow\RR$ the measurable function
$$Qf(\omega)=\int_\mathbf{\Omega'} f(\widetilde{\omega})Q(d\widetilde{\omega}|\omega)\quad\hbox{for $\omega\in\mathbf{\Omega}$}.$$
For a measure $\mu\in\bcal{M}^{+}(\mathbf{\Omega})$, we denote by $\mu Q$ the finite measure  on $(\widetilde{\mathbf{\Omega}},\widetilde{\mathcal{F}})$ given by
$$B\mapsto \mu Q\,(B)= \int_{\mathbf{\Omega}} Q(B|\omega) \mu(d\omega)\quad\hbox{for $B\in\widetilde{\mathcal{F}}$}.$$
The product of the $\sigma$-algebras $\mathcal{F}$ and $\widetilde{\mathcal{F}}$ is denoted by $\mathcal{F}\otimes\widetilde{\mathcal{F}}$ and consists of the $\sigma$-algebra
generated by the measurable rectangles, that is, the sets of the form $\Gamma\times\widetilde{\Gamma}$ for $\Gamma\in\mathcal{F}$ and
$\widetilde{\Gamma}\in\widetilde{\mathcal{F}}$.
We denote by $\mu\otimes Q$ the unique  finite measure on the product space $(\mathbf{\Omega}\times\widetilde{\mathbf{\Omega}},\mathcal{F}\otimes\widetilde{\mathcal{F}})$ satisfying 
$$(\mu\otimes Q)( \Gamma\times\widetilde{\Gamma})= \int_{\Gamma} Q(\widetilde{\Gamma}|\omega)\mu(d\omega)\quad\hbox{for $\Gamma\in\mathcal{F}$ and
$\widetilde{\Gamma}\in\widetilde{\mathcal{F}}$,}$$
see Proposition III-2-1 in \cite{neveu70} for a proof of the existence and uniqueness of such measure.
Let $(\widebar{\mathbf{\Omega}},\widebar{\mathcal{F}})$ be a third measurable space and $R$ be a stochastic kernel on $\widebar{\mathbf{\Omega}}$
given $\widetilde{\mathbf{\Omega}}$. Then we will denote by $QR$ the stochastic kernel on $\widebar{\mathbf{\Omega}}$ given $\mathbf{\Omega}$ defined as
$$QR(\Gamma |\omega)= \int_{\widetilde{\mathbf{\Omega}}} R(\Gamma | \tilde{\omega}) Q(d\tilde{\omega} | \omega)
\quad\hbox{for $\Gamma\in\widebar{\mathcal{F}}$ and $\omega\in\mathcal{F}$}.$$
Given $\mu\in\bcal{M}(\mathbf{\Omega}\times\widetilde{\mathbf{\Omega}})$,  the  marginal measures are $\mu^{\mathbf{\Omega}}\in\bcal{M}(\mathbf{\Omega})$ and
$\mu^{\widetilde{\mathbf{\Omega}}}\in\bcal{M}(\widetilde{\mathbf{\Omega}})$ defined by 
$\mu^{\mathbf{\Omega}}(\cdot)=\mu(\cdot\times \widetilde{\mathbf{\Omega}})$ and
$\mu^{\widetilde{\mathbf{\Omega}}}(\cdot)=\mu(\mathbf{\Omega}\times\cdot)$. If $\pi$ is a kernel on $\widetilde{\mathbf{\Omega}}\times\widebar{\mathbf{\Omega}}$ given $\mathbf{\Omega}$ the  marginal kernels are $\pi^{\widetilde{\mathbf{\Omega}}}$
and $\pi^{\widebar{\mathbf{\Omega}}}$, respectively defined by 
$\pi^{\widetilde{\mathbf{\Omega}}}=\pi(\cdot\times \widebar{\mathbf{\Omega}}|\omega)$ and
$\pi^{\widebar{\mathbf{\Omega}}}=\pi(\widetilde{\mathbf{\Omega}}\times\cdot |\omega)$ for $\omega\in\mathbf{\Omega}$.

Given metric spaces $\SS$ and $\SS'$ (endowed with their respective Borel $\sigma$-algebras),
we say that $f:\mathbf{\Omega}\times\mathbf{S}\rightarrow\mathbf{S}'$ is a \textit{Carath\'eodory function} if $f(\cdot,s)$ is measurable on~$\mathbf{\Omega}$ for every $s\in\mathbf{S}$ and $f(\omega,\cdot)$ is continuous on $\mathbf{S}$ for every $\omega\in \mathbf{\Omega}$. The family of the so-defined Carath\'eodory functions is denoted by $\car(\mathbf{\Omega}\times\mathbf{S},\mathbf{S}')$. The family of Carath\'eodory functions which, in addition, are bounded is denoted by $\car_b(\mathbf{\Omega}\times\mathbf{S},\mathbf{S}')$. 
When the metric space~$\mathbf{S}$ is separable then any $f\in\car(\mathbf{\Omega}\times\mathbf{S},\mathbf{S}')$ is a jointly measurable function on $(\mathbf{\Omega}\times\mathbf{S},\mathcal{F}\otimes\bfrak{B}(\mathbf{S}))$; see \cite[Lemma 4.51]{aliprantis06}.

If $\mathbf{S}$ is a Polish space (a complete and separable metric space), on $\bcal{M}(\mathbf{\Omega}\times\mathbf{S})$ we will consider the  $ws$-topology (weak-strong topology) which is the coarsest topology for which the mappings
$$
\mu\mapsto \int_{\mathbf{\Omega}\times\mathbf{S}} f(\omega,s)\mu(d\omega,ds)
$$
for $f\in\car_b(\mathbf{\Omega}\times\mathbf{S},\RR)$ are continuous.
There are other equivalent definitions of this topology as discussed, for instance, in \cite[Section 3.3]{florescu12}.

\bigskip

The next disintegration lemma will be useful in the forthcoming (see Theorem 1 in \cite{valadier73}).
\begin{lemma}[Disintegration lemma]
\label{lemma-disintegration} 
Let $(\mathbf{\Omega},\mathcal{F})$ be a measurable space and let $\mathbf{S}$ be a Polish space. Let $\Phi:\mathbf{\Omega}\to 2^{\mathbf{S}}$ be a weakly measurable correspondence with nonempty closed values, and let~$\mathbf{K}$ be the graph of the correspondence. For every  $\mu\in\bcal{M}^{+}(\mathbf{\Omega}\times\mathbf{S})$ 
such that $\mu(\mathbf{K}^{c})=0$
there exists a stochastic kernel $Q$ on $\mathbf{S}$ given $\mathbf{\Omega}$ such that 
\begin{equation}\label{eq-dudley-product}
\mu= \mu^{\mathbf{\Omega}}\otimes Q 
\end{equation}
and such that
$Q(\Phi(\omega)|\omega)=1$ for each $\omega\in\mathbf{\Omega}$. Moreover, $Q$ is unique $\mu^{\mathbf{\Omega}}$-almost surely, meaning that if $Q$ and $Q'$ are two stochastic kernels that satisfy \eqref{eq-dudley-product} then for all $\omega$ in a set of $\mu^{\mathbf{\Omega}}$-probability one, the probability measures $Q(\cdot|\omega)$ and $Q'(\cdot|\omega)$ coincide.
\end{lemma}

%%%%%%%%%%%%%%%%%%%%%%%%%%%%%%%%%%%%%%%
\section{The absorbing control model}
\label{Sec-2}
%%%%%%%%%%%%%%%%%%%%%%%%%%%%%%%%%%%%%%%
The main goal of this section is to introduce the parameters defining the model with a brief presentation of the construction of the controlled process.
We also describe the notions of  absorbing and uniformly absorbing MDP, and give the definition of an occupation measure, providing a first elementary property (see Lemma \ref{cor-bounded}).

\subsection{The control model.}
We consider a stationary Markov controlled  process $(\mathbf{X},\mathbf{A},\{\mathbf{A}(x): x\in \mathbf{X}\},Q,\eta)$
consisting of:
\begin{itemize}
\item A measurable state space $\XX$ endowed with a $\sigma$-algebra $\bfrak{X}$.
\item A Borel space $\mathbf{A}$, representing the action space.
\item A family of nonempty measurable sets $\AA(x)\subseteq \AA$ for $x\in\XX$. The set $\AA(x)$ gives the available actions in state~$x$.
Let $\KK=\{(x,a)\in\XX\times\AA: a\in \AA(x)\}$ be the family of feasible state-action pairs. We assume that $\KK$ is a measurable subset of
$\XX\times\AA$ endowed with the $\sigma$-algebra $\bfrak{X}\otimes\bfrak{B}(\AA)$.
We also assume that $\KK$ contains the graph of a measurable selector (i.e., the graph of a measurable $f:\XX\rightarrow\AA$ such that $f(x)\in\AA(x)$ for every $x\in\XX$).
\item A stochastic kernel $Q$ on $\mathbf{X}$ given $\XX\times\AA$, which stands for the transition probability function. 
\item An initial distribution given by $\eta\in\bcal{P}(\XX)$.
\end{itemize}
Additionally, we assume that a measurable set  $\Delta\in\bfrak{X}$ is given. As we shall see later in Definition~\ref{absorbing}, the set $\Delta$ will play the role of the \emph{absorbing set}. 
The so-defined control model is denoted by $\bfrak{M}(\eta,\Delta)$ where we make explicit the dependence on the initial distribution and the absorbing set. 

\bigskip

The space of admissible histories of the controlled process up to time $n\in\NN$ is denoted by $\HH_{n}$. It is defined recursively by $$\HH_0=\XX\quad\hbox{and}\quad
 \HH_n=\HH_{n-1}\times\AA\times\XX\quad\hbox{for $n\ge1$},$$ all endowed with their corresponding product $\sigma$-algebras.
 A control policy $\pi$ is a sequence $\{\pi_n\}_{n\ge0}$ of stochastic kernels on $\AA$ given $\HH_n$, denoted by $\pi_n(da|h_n)$, such that
$$\pi_n(\AA(x_n)|h_n)=1\quad\hbox{for each $n\ge0$ and $h_n=(x_0,a_0,\ldots,x_n)\in\HH_n$}.$$  
The set of all policies is denoted by $\mathbf{\Pi}$. 

Let us denote by $\mathbf{M}$ the set of stochastic kernels $\sigma$ on $\AA$ given $\XX$ satisfying $\sigma(\AA(x)|x)=1$
for any $x\in \XX$.
A policy $\pi=\{\pi_{n}\}_{n\in\NN}\in \mathbf{\Pi}$ is called a \emph{stationary randomized policy} if there exists $\sigma\in\mathbf{M}$ satisfying
$$\pi_{n}(\cdot |h_{n})=\sigma(\cdot |x_{n})\quad\hbox{ for any $h_{n}=(x_{0},a_{0},\ldots,x_{n})\in \mathbf{H}_{n}$ and $n\in \NN$}.$$
In such a case, we will write $\sigma$ instead of $\pi$ to emphasize that the corresponding stationary randomized policy $\pi$ is generated by $\sigma$.
Since the set of stationary randomized policies can be identified with $\mathbf{M}$, 
will refer to it as $\mathbf{M}$ in a slight abuse of notation.
%We denote by $\mathbf{\Pi}^{s}$ the set of all stationary randomized policies.
%A policy $\pi=\{\pi_{n}\}_{n\in\NN}\in \mathbf{\Pi}$ is called a \emph{Markovian randomized policy} if there exists a sequence $\{\varphi_{n}\}_{n\in\NN} $ in $\mathbf{M}$ satisfying
%$$\pi_{n}(\cdot |h_{n})=\varphi_{n}(\cdot |x_{n})\quad\hbox{ for any $h_{n}=(x_{0},a_{0},\ldots,x_{n})\in \mathbf{H}_{n}$ and $n\in \NN$}.$$
%Let $\mathbf{\Pi}^{m}$ be the set of all Markovian randomized policies.
%We have  $\mathbf{\Pi}^s\subseteq\mathbf{\Pi}$.

For each $\sigma\in\mathbf{M}$, we denote by $Q_\sigma$ the stochastic kernel on $\XX$ given $\XX$ defined by
\begin{equation*}
Q_{\sigma}(D|x)= \int_{\AA} Q(D|x,a)\sigma(da|x)
\quad\hbox{for $x\in\XX$ and $D\in\bfrak{X}$}.
\end{equation*}
The compositions of $Q_\sigma$ with itself are denoted by $Q^t_{\sigma}$ for any $t\ge0$, with the convention that $Q^0_\sigma(\cdot|x)$ is the Dirac probability measure concentrated at $x$. 

\bigskip

The canonical space of all possible sample paths of the state-action process is 
$\mathbf{\Omega}=(\XX\times\AA)^{\infty}$ endowed with the product $\sigma$-algebra $\bcal{F}$. 
The coordinate projection functions from $\mathbf{\Omega}$ to the state space $\XX$, the action space $\AA$, and $\HH_n$ for $n\ge0$ are respectively denoted   by $X_{n}$, $A_{n}$, and $H_n$. 
We will refer  to $\{X_{n}\}_{n\in \NN}$ as to state process and $\{A_{n}\}_{n\in \NN}$ as the action  process.
It is a well known result that for every policy $\pi \in \mathbf{\Pi}$ and any initial probability measure $\lambda$ on $(\mathbf{X},\bfrak{X})$ there exists a unique probability
measure $\mathbb{P}_{\lambda,\pi}$ on $(\mathbf{\Omega},\bcal{F})$
such that $ \mathbb{P}_{\lambda,\pi} (\KK^{\infty})=1$ and such that for every $n\in\NN$, $\Gamma\in \bfrak{X}$, and  $ \Lambda\in \bfrak{B}(\mathbf{A})$ we have 
$\mathbb{P}_{\lambda,\pi}(X_{0}\in \Gamma)=\lambda(\Gamma)$, 
$$\mathbb{P}_{\lambda,\pi}(X_{n+1}\in \Gamma\mid H_{n},A_{n}) = Q(\Gamma\mid X_{n},A_{n}) \quad\hbox{and}\quad
\mathbb{P}_{\lambda,\pi}(A_{n}\in \Lambda\mid H_{n})= \pi_{n}(\Lambda\mid H_{n})$$
with
$\mathbb{P}_{\lambda,\pi}$-probability one.
The expectation with respect to  $\mathbb{P}_{\lambda,\pi}$ is denoted by $\mathbb{E}_{\lambda,\pi}$.

\begin{definition}
The hitting time $T_\Delta$ of the set $\Delta$ is given by  $T_\Delta:\mathbf{\Omega}\rightarrow\NN\cup\{\infty\}$ defined as
$$T_\Delta(x_0,a_0,x_1,a_1,\ldots)=\min\{n\ge0: x_n\in\Delta\},$$
where the $\min$ over the empty set is defined as $+\infty$. \end{definition}

Next we propose the definition of an absorbing control model.

\begin{definition}
\label{absorbing}
Given an initial distribution $\lambda\in\bcal{P}(\XX)$ and $\Delta\in\bfrak{X}$, we say that the control model  
$\bfrak{M}(\lambda,\Delta)$ is absorbing  if the conditions (a) and (b) below are satisfied, and we say that 
it is uniformly absorbing if, additionally, condition (c) holds.
\begin{enumerate}[label=(\alph*).]
\item 
For every $(x,a)\in\Delta\times\AA$ we have $Q(\Delta|x,a)=1$. 
\item For any $\pi\in\mathbf{\Pi}$ the expected hitting time $\mathbb{E}_{\lambda,\pi}[T_\Delta]$ is finite.
\item The following limit holds: $$\lim_{n\rightarrow\infty} \sup_{\pi\in\mathbf{M}} \sum_{t=n}^\infty\mathbb{P}_{\lambda,\pi}\{T_\Delta>t\}=0.$$
\end{enumerate}
\end{definition}

We define now the occupation measures of an absorbing control model $\bfrak{M}(\lambda,\Delta)$. 

\begin{definition}\label{def-occupation-measure}
The  occupation measure   $\mu_{\lambda,\pi}$  of the control policy $\pi\in\mathbf{\Pi}$ is 
\begin{eqnarray*}
\mu_{\lambda,\pi}(\Gamma) &=& \mathbb{E}_{\lambda,\pi}\Big[ \sum_{t=0}^\infty \mathbf{I}_{\{T_\Delta>t\}}\cdot\mathbf{I}_{\{(X_t,A_{t})\in \Gamma\}}\Big] \quad\hbox{for $\Gamma\in\bfrak{X}\otimes\mathfrak{B}(\AA)$.}
\end{eqnarray*}
\end{definition}

We note that the occupation measure $\mu_{\lambda,\pi}$ takes into account the state-action process up to time $T_\Delta$ and it does not count the time spent in $\Delta$: indeed, $\mu_{\lambda,\pi}(\Delta\times\AA)=0$.  
We will consider the following sets of occupation measures 
$\bcal{O}_\lambda=\{\mu_{\lambda,\pi}:\pi\in\mathbf{\Pi}\}$ %, $\bcal{O}^m_\lambda=\{\mu_{\lambda,\pi}:\pi\in\mathbf{\Pi}^m\}$
and also $\bcal{O}^s_\lambda=\{\mu_{\lambda,\pi}:\pi\in\mathbf{M}\}$.

\begin{lemma}\label{cor-bounded} 
For an absorbing control model $\bfrak{M}(\lambda,\Delta)$, the set $\bcal{O}_\lambda$ is included in $ \MM^{+}(\XX\times\AA)$ and it is bounded.
\end{lemma}
\textbf{Proof.}
Note that the total mass of the occupation measure $\mu_{\lambda,\pi}$ is 
$$\mu_{\lambda,\pi}(\XX\times\AA)=\mathbb{E}_{\lambda,\pi}\Big[\sum_{t=0}^\infty \mathbf{I}_{\{T_\Delta>t\}}\Big]
=\sum_{t=0}^\infty \mathbb{P}_{\lambda,\pi}\{T_\Delta>t\}=\mathbb{E}_{\lambda,\pi}[T_\Delta],
$$
which is finite as a consequence of item (b) in Definition \ref{absorbing}. We also have $\sup_{\pi\in\mathbf{\Pi}}\mathbb{E}_{\lambda,\pi}[T_\Delta]<\infty$ according to \cite[Sections 4.4 and 5.5]{dynkin79} for the special case of a Borel state space or Proposition 2.4(i) in \cite{dufour23} for the general case of a measurable state space. This shows that the set $\bcal{O}_\lambda$ is bounded.
\hfill$\Box$

\section{The characteristic equations}
\label{Sec-3}
In this section, we introduce the notion of phantom measure and show that, for an absorbing MDP, phantom measures may exist.
In order to avoid this pathological phenomenon, we establish our first main result (see Theorem \ref{prop-characterizations}) that provides several necessary and sufficient conditions to ensure that there are no phantom measures.

\begin{definition} Given the absorbing control model $\bfrak{M}(\eta,\Delta)$,
we say  $\mu\in\MM^{+}(\XX\times\AA)$ is a solution of the characteristic equations when 
\begin{equation}\label{eq-linear-equation} 
\mu(\KK^c)=0\quad\hbox{and}\quad \mu^\XX=(\eta+\mu Q)\mathbb{I}_{\Delta^c}.
\end{equation}
We denote by $\bcal{C}_\eta$ the family of all solutions of the characteristic equations.
\end{definition}

\begin{lemma}
\label{vanishing}
Suppose that  the control model $\bfrak{M}(\eta,\Delta)$ is absorbing. Given any $\pi\in\mathbf{\Pi}$ and $\sigma\in\mathbf{M}$, we have
$$\mathbb{P}_{x,\sigma}\{T_\Delta=\infty\}=0 \text{ for } \mu_{\eta,\pi}^\XX \text{-almost all } x\in\XX.$$
\end{lemma}
\textbf{Proof.}
For notational convenience, let us introduce the functions $h_t$ and $h$ on $\XX$ taking values in $[0,1]$ as
$h_t(x)=\mathbb{P}_{x,\sigma}\{T_{\Delta} > t\}$ and $h(x)=\mathbb{P}_{x,\sigma}\{T_\Delta=\infty\}$ respectively.
Note that $h_t$ vanishes on $\Delta$ and also that $h_t$  is a decreasing sequence of measurable functions which converges pointwise to 
$h(x)=\mathbb{P}_{x,\sigma}\{T_\Delta=\infty\}$.   
To prove the result, we proceed by contradiction and suppose that $\int h d\mu^\XX_{\eta,\pi}>0$. By definition of the occupation measure
$\mu_{\eta,\pi}$ and since $h$ vanishes on~$\Delta$, this implies that there exists some  $s\ge0$ with $\mathbb{E}_{\eta,\pi}[h(X_s)]>0$ and, in particular, since $h_t\downarrow h$,
\begin{equation}\label{eq-lemma-series-1}
\sum_{t=0}^\infty \mathbb{E}_{\eta,\pi}[h_t(X_s)]=\infty.
\end{equation}
Define the strategy $\gamma\in\mathbf{\Pi}$ as follows: $\gamma_k=\pi_k$ for $0\le k\le s-1$ and
 $\gamma_{k}=\sigma$ for every  $k\geq s$ (in case that $s=0$, this definition reduces to  $\gamma=\sigma$).
Consequently, the distribution of $X_{s}$ is the same under $\mathbb{P}_{\eta,\pi}$ and $\mathbb{P}_{\eta,\gamma}$, and we have 
$$\mathbb{E}_{\eta,\pi}[h_t(X_s)]=\mathbb{E}_{\eta,\gamma}[h_t(X_s)]=\mathbb{P}_{\eta,\gamma}\{T_\Delta>t+s\}.$$
Recalling \eqref{eq-lemma-series-1}, this implies that the series
$$\sum_{t=0}^\infty\mathbb{P}_{\eta,\gamma}\{T_\Delta>t+s\}=\infty,$$ in contradiction with $\mathbb{E}_{\eta,\gamma}[T_\Delta]$ being finite. 
This gives the result.
$\hfill$ $\Box$
%\\[10pt]\indent

\begin{proposition}\label{prop-linear-equation-occupation} 
Suppose that  the control model $\bfrak{M}(\eta,\Delta)$ is absorbing.
\begin{itemize}
\item[(i).] Given any $\pi\in\mathbf{\Pi}$, its occupation measure $\mu_{\eta,\pi}$ satisfies the characteristic equations. Therefore, we have the inclusion $\bcal{O}_\eta\subseteq\bcal{C}_\eta$.
\item[(ii).] Given any $\pi\in\mathbf{\Pi}$, there exists $\sigma\in\mathbf{M}$ such that $\mu_{\eta,\pi}=\mu_{\eta,\sigma}$. Hence, $\bcal{O}^s_\eta=\bcal{O}_\eta$.
\end{itemize}
\end{proposition}
\textbf{Proof.} (i). We have already shown that $\mu_{\eta,\pi}$ is a finite measure on $\XX\times\AA$.
To prove the stated result, note that for any $B\in\bfrak{X}$ we have
\begin{align*}
\mu^\XX(B) = \sum_{t=0}^\infty \mathbb{P}_{\eta,\pi}\{T_\Delta>t,X_t\in B\} 
=  \eta(B-\Delta)+\sum_{t=1}^\infty \mathbb{E}_{\eta,\pi}\big[  \mathbb{P}_{\eta,\pi}\{T_\Delta>t,X_t\in B\mid H_{t-1},A_{t-1}\}\big].
\end{align*}
Observe now that for each $t\ge1$  the conditional probability within brackets vanishes on the set $\{T_\Delta\le t-1\}$, and so
\begin{eqnarray*}
\mu^\XX(B) 
&= & \eta(B-\Delta)+\sum_{t=1}^\infty \mathbb{E}_{\eta,\pi}\big[  Q(B-\Delta\mid X_{t-1},A_{t-1})
\cdot\mathbf{I}_{\{T_\Delta>t-1\}}\big]\\
&=&  \eta(B-\Delta)+\int_{\XX\times\AA}Q(B-\Delta |x,a)\mu(dx,da),
\end{eqnarray*}
which can be equivalently written precisely as $\mu^\XX=(\eta+\mu Q)\mathbb{I}_{\Delta^c}$.
By construction of the state-action process, it is clear that $\mu(\KK^c)=0$. This shows that $\mu_{\eta,\pi}$ indeed verifies the characteristic equations. 
\\[5pt]\noindent (ii). We will show that for any $\pi\in\mathbf{\Pi}$ there is some $\sigma\in\mathbf{M}$ with $\mu_{\eta,\pi}=\mu_{\eta,\sigma}$.
Since $\mu_{\eta,\pi}$ is in $\MM^+(\XX\times\AA)$ with $\mu_{\eta,\pi}(\KK^c)=0$, the measure $\mu_{\eta,\pi}$ can be disintegrated as 
\begin{equation}\label{eq-lemma-series-2}
\mu_{\eta,\pi}=\mu_{\eta,\pi}^\XX\otimes \sigma\quad\hbox{for some $\sigma\in\mathbf{M}$}.
\end{equation}
It follows that the characteristic equation can be written as
$$
\mu_{\eta,\pi}^\XX= (\eta+  \mu_{\eta,\pi} Q )\mathbb{I}_{\Delta^c} = (\eta+\mu_{\eta,\pi}^\XX Q_\sigma)\mathbb{I}_{\Delta^c}.
$$
Iterating this equation  we obtain that for any $t\ge0$ 
\begin{align}
\mu_{\eta,\pi}^\XX=\eta \mathbb{I}_{\Delta^c} \sum_{k=0}^{t-1} (Q_{\sigma}\mathbb{I}_{\Delta^c})^k + \mu_{\eta,\pi}^\XX ( Q_{\sigma}\mathbb{I}_{\Delta^c} )^{t}.
\label{expression-marginal-1}
\end{align}
We have, for any $x\in\XX$ and $t\ge0$,
$$( Q_{\sigma}\mathbb{I}_{\Delta^c} )^{t}(\XX |x)=Q_\sigma^t(\Delta^c|x)=\mathbb{P}_{x,\sigma}\{T_{\Delta} > t\}$$ 
and $\mathbb{P}_{x,\sigma}\{T_{\Delta} > t\} \downarrow \mathbb{P}_{x,\sigma}\{T_\Delta=\infty\}$ as $t\rightarrow\infty$ for any $x\in\XX$.
Applying Lemma \ref{vanishing}, we obtain that $ \lim_{t\rightarrow\infty}\mu_{\eta,\pi}^\XX ( Q_{\sigma}\mathbb{I}_{\Delta^c} )^{t}(\XX)=0$.
Taking the limit as $t\rightarrow\infty$  in \eqref{expression-marginal-1} we obtain 
$$
\mu_{\eta,\pi}^\XX=\eta \mathbb{I}_{\Delta^c} \sum_{k=0}^{\infty} (Q_{\sigma}\mathbb{I}_{\Delta^c})^k =\mu_{\eta,\sigma}^\XX.
$$
Once we have shown that the $\XX$-marginals of the occupation measures of $\mu_{\eta,\pi}$ and $\mu_{\eta,\sigma}$ coincide, we conclude the result from \eqref{eq-lemma-series-2}.
$\hfill$ $\Box$
\begin{definition}
A measure $\mu\in\MM^{+}(\XX\times\AA)$ is called a phantom solution to the characteristic equation or simply a phantom measure if it satisfies the characteristic equations \eqref{eq-linear-equation} but it is not given by an the occupation measure of any policy in~$\mathbf{\Pi}$, that is, $\mu\in \bcal{C}_\eta$ but $\mu\notin \bcal{O}_\eta$.
\end{definition}
The following simple example illustrates the existence of phantom measures and, in particular, the existence of such measures  shows that \cite[Lemma 4.3]{feinberg12} is not stated in a precise way.
\begin{example}\label{example-phantom}
Let us consider the model defined by $\mathbf{X}=\{0,1,2\}$, $\mathbf{A}=\{a\}$, $\Delta=\{0\}$, $\eta=\delta_{\{1\}}$, and 
$$Q(\{0\}|0,a)=Q(\{0\}|1,a)=Q(\{2\}|2,a)=1.$$
A measure $\nu$ on $\XX$ is a solution of \eqref{eq-linear-equation} if and only if 
\[
\nu(\{0\})=0,\  \nu(\{1\})=1,\ \hbox{and}\ \nu(\{2\})=K,\]
where $K$ is an arbitrary nonnegative real number.
The unique occupation measure of this model corresponds to $K=0$ while, for $K>0$,  we obtain phantom measures.
It will be shown in Theorem~\ref{Eq-Characteristic++}
that the occupation measures are characterized by \eqref{eq-linear-equation} and an additional condition.
This is the first notable difference with the discounted model, where the occupation measures are essentially characterized by a single equation of type
\eqref{eq-linear-equation}, that is, $\mu(\KK^c)=0$ and $\mu^\XX=\eta+  \alpha\mu Q $ (where $\alpha$ is the discount factor), thus
excluding the existence of phantom measures.
\end{example}
\begin{definition}\label{def-invariant} 
Let $\bfrak{M}(\eta,\Delta)$ be an absorbing control model.
Given a measure $\vartheta\in\MM^{+}(\XX\times\AA)$ we say that $\vartheta$ is $Q\mathbb{I}_{\Delta^c}$-invariant when 
$\vartheta(\KK^c)=0$ and $ \vartheta^\XX=\vartheta Q\mathbb{I}_{\Delta^c}$.
\end{definition}

Let us make some comments on this definition. Note that we are not excluding that $\vartheta$ is the null measure.
We call such a measure invariant because, by disintegration, there exists $\sigma\in\mathbf{M}$ satisfying $\vartheta^\XX=\vartheta^\XX Q_\sigma\mathbb{I}_{\Delta^c}$ which can be written as 
$$\vartheta^\XX(\Delta)=0\quad\hbox{and}\quad \vartheta^\XX(B)=\int_\XX Q_\sigma(B|x)\vartheta^\XX(dx)$$
for measurable $B\subseteq\Delta^c$, and so $\vartheta^\XX$ is an invariant measure for the substochastic kernel 
$Q_\sigma$ on~$\Delta^c$.

\begin{theorem}\label{prop-characterizations} 
Let $\bfrak{M}(\eta,\Delta)$ be absorbing. 
\begin{itemize}
\item[(i).]
A  measure  is in $\bcal{C}_\eta$ if and only if  it can be decomposed as the sum of an occupation measure $\mu_{\eta,\sigma}$ in~$\bcal{O}_\eta^s$
and a $Q\mathbb{I}_{\Delta^c}$-invariant measure $\nu\otimes\sigma$ with $\sigma\in\mathbf{M}$ and $\nu\in\MM^{+}(\XX)$.
\item[(ii).] The following statements are equivalent.
\begin{itemize}
\item[(a).] The unique $Q\mathbb{I}_{\Delta^c}$-invariant measure is the null measure on $\XX\times\AA$. 
\item[(b).] $\bcal{C}_\eta=\bcal{O}_\eta$, i.e., there are no phantom measures.
\item[(c).] The set $\bcal{C}_\eta$ is bounded.
\end{itemize}
\end{itemize}
\end{theorem}
\textbf{Proof.} (i).
Suppose that $\mu\in\bcal{M}^+(\XX\times\AA)$ is a solution of the characteristic equations. Proceeding as in the proof of Proposition \ref{prop-linear-equation-occupation}, we derive the existence of $\sigma\in\mathbf{M}$  such that
$\mu=\mu^\XX\otimes\sigma$ and 
\begin{align}
\mu^\XX=\eta \mathbb{I}_{\Delta^c} \sum_{k=0}^{t-1} (Q_{\sigma}\mathbb{I}_{\Delta^c})^k + \mu^\XX ( Q_{\sigma}\mathbb{I}_{\Delta^c} )^{t}
\end{align}
for any $t\ge1$. It follows that by taking the limit as $t\rightarrow\infty$ that 
$$\nu=\lim_{t\rightarrow\infty}\mu^\XX ( Q_{\sigma}\mathbb{I}_{\Delta^c} )^{t}=\mu^\XX-\mu_{\eta,\sigma}^\XX$$
defines a measure in $\MM^{+}(\XX\times\AA)$ satisfying $\mu^\XX=\mu_{\eta,\sigma}^\XX+\nu$.
Recalling that $$\mu^\XX=\big(\eta+\mu^\XX Q_{\sigma}\big) \mathbb{I}_{\Delta^c}\quad\hbox{and}\quad
\mu_{\eta,\sigma}^\XX=\big(\eta+\mu_{\eta,\sigma}^\XX Q_{\sigma}\big) \mathbb{I}_{\Delta^c},$$
we get that
$\nu=\nu Q_{\sigma}\mathbb{I}_{\Delta^c}$. The measure $\vartheta=\nu\otimes\sigma$ satisfies the conditions in the statement of this proposition: $\mu=\mu_{\eta,\sigma}+\vartheta$. 

Conversely, it is straightforward to check that the sum of a $Q\mathbb{I}_{\Delta^c}$-invariant measure and an occupation measure (which satisfies the characteristic equations)  satisfies itself the characteristic equations and, hence, belongs to
$\bcal{C}_\eta$. In fact, the sum of any measure in $\bcal{C}_\eta$ and a $Q\mathbb{I}_{\Delta^c}$-invariant measure lies in $\bcal{C}_\eta$.
\\[5pt]
\noindent
(ii). The implication $ (a)\Rightarrow(b)$  follows directly from item (i), while $(b)\Rightarrow(c)$ is derived from Lemma \ref{cor-bounded}.
We prove $(c)\Rightarrow(a)$ by contradiction and, hence, we suppose that there exists a non-null $Q\mathbb{I}_{\Delta^c}$-invariant measure $\vartheta$.
For any $K>0$ we have that $K\vartheta$ is also a
$Q\mathbb{I}_{\Delta^c}$-invariant measure and so for any $\mu\in\bcal{C}_\eta$ and every $K>0$ we have $\mu+K\vartheta\in\bcal{C}_\eta$, which is not compatible with $\bcal{C}_\eta$ being bounded. 
 \hfill $\Box$
\\[10pt]\indent
Our next result gives some insight on the behavior of phantom measures. In particular, the first part of the result establishes, loosely speaking, that the state process $\{X_t\}$ under any policy $\pi\in\mathbf{\Pi}$ never visits the support of a non-null $Q\mathbb{I}_{\Delta^c}$-invariant measure.
At this point, recall that a $Q\mathbb{I}_{\Delta^c}$-invariant measure is in $\MM^{+}(\XX\times\AA)$ by definition, and so it is a finite measure.
\begin{proposition}
\label{Prop-invariant}
Suppose that $\bfrak{M}(\eta,\Delta)$ is absorbing. 
\begin{itemize}
\item[(i).] Let $\vartheta$ be a non-null $Q\mathbb{I}_{\Delta^c}$-invariant measure, we have $\mu_{\eta,\pi}^\XX\perp \vartheta^\XX$ for every $\pi\in\mathbf{\Pi}$.
\item[(ii).] For any $\mu\in\bcal{C}_{\eta}$, there exist a unique $\mu_{\eta,\sigma}\in\bcal{O}_\eta^s$ and a unique $Q\mathbb{I}_{\Delta^c}$-invariant measure $\vartheta$ such that
$$\mu=\mu_{\eta,\sigma}+\vartheta.$$
\end{itemize}
\end{proposition}
\textbf{Proof.} 
(i). By disintegration of $\vartheta$ there exists $\sigma\in\mathbf{M}$ with
$\vartheta^\XX=\vartheta^\XX Q_\sigma\mathbb{I}_{\Delta^c}$. Iterating this equation we obtain $\vartheta^\XX=\vartheta^\XX Q^t_\sigma\mathbb{I}_{\Delta^c}$ for every $t\ge1$ and taking the limit as $t\rightarrow\infty$ yields
$$\vartheta^\XX(\XX)=\vartheta^\XX(\Delta^c)=\int_\XX \mathbb{P}_{x,\sigma}\{T_\Delta=\infty\} \vartheta^\XX(dx)$$
and so the set $B=\{x\in\XX:  \mathbb{P}_{x,\sigma}\{T_\Delta=\infty\}<1\}$ satisfies $\vartheta^\XX(B)=0$.
By the Lebesgue decomposition theorem there exist two finite measures $\vartheta_1^\XX,\vartheta_2^\XX\in\MM_+(\XX)$ with
$$\vartheta^\XX= \vartheta_1^\XX+\vartheta_2^\XX\quad\hbox{with $\vartheta_1^\XX\ll \mu_{\eta,\pi}^\XX$ and
$\vartheta_2^\XX\perp \mu_{\eta,\pi}^\XX$}.$$
Therefore,
\begin{eqnarray*}
\vartheta^\XX(\XX) &=& \int_\XX \mathbb{P}_{x,\sigma}\{T_\Delta=\infty\} \vartheta^\XX_1(dx)+\int_\XX \mathbb{P}_{x,\sigma}\{T_\Delta=\infty\} \vartheta^\XX_2(dx)\\
&=&\int_\XX \mathbb{P}_{x,\sigma}\{T_\Delta=\infty\}\cdot \frac{d\vartheta^\XX_1}{d\mu^\XX_{\eta,\pi}}\ \mu^\XX_{\eta,\pi}(dx)+\int_\XX \mathbb{P}_{x,\sigma}\{T_\Delta=\infty\} \vartheta^\XX_2(dx).
\end{eqnarray*}
Applying Lemma \ref{vanishing}, we obtain that $\mathbb{P}_{x,\sigma}\{T_\Delta=\infty\}$ vanishes $\mu^\XX_{\eta,\pi}$-a.s. and so 
$$\vartheta^\XX(\XX) = \int_\XX \mathbb{P}_{x,\sigma}\{T_\Delta=\infty\} \vartheta^\XX_2(dx).$$
Observe that $\vartheta^\XX_2(B)=0$ and so
$$
\vartheta^\XX(\XX) =  \int_{B^c} \mathbb{P}_{x,\sigma}\{T_\Delta=\infty\} \vartheta^\XX_2(dx) 
= \vartheta^\XX_2(B^c)\ =\ \vartheta^\XX_2(\XX).
$$
This shows that the measure $\vartheta^\XX_1$ is null and, therefore, $\vartheta^\XX=\vartheta^\XX_2$, showing the result.
\\[5pt]
\noindent
(ii). Consider $\mu\in\bcal{C}_\eta$. Suppose that 
\begin{align}
\mu=\mu_{\eta,\sigma}+\vartheta_{\sigma}=\mu_{\eta,\gamma}+\vartheta_{\gamma}
\end{align}
for $\mu_{\eta,\sigma}\in\bcal{O}_\eta^s$ and $\mu_{\eta,\gamma}\in\bcal{O}_\eta^s$, and where $\vartheta_{\sigma}$ and $\vartheta_{\gamma}$ are$Q\mathbb{I}_{\Delta^c}$-invariant measures.
By using standard arguments used for the special case of a Borel state space,
it can be shown that the set of strategic measures $ \big\{\mathbb{P}_{\eta,\pi}\}_{\pi\in\mathbf{\Pi}}$ is convex for the  general case of a measurable state space. 
This implies the convexity of $\bcal{O}_\eta$. Moreover, the set of $Q\mathbb{I}_{\Delta^c}$-invariant measures is convex by definition.
Therefore, by observing that $\mu=\frac{1}{2}(\mu_{\eta,\sigma}+\mu_{\eta,\gamma}) + \frac{1}{2}(\vartheta_{\sigma}+\vartheta_{\gamma})$ and from item (i), there exists a set
$\Lambda\in\bfrak{X}$ such that $\mu_{\eta,\sigma}(\Lambda^{c}\times\AA) = \mu_{\eta,\gamma}(\Lambda^{c}\times\AA)=0$ and
$\vartheta_{\sigma}(\Lambda\times\AA)=\vartheta_{\gamma}(\Lambda\times\AA)=0$. 
Therefore, $\mu_{\eta,\sigma} = \mu_{\eta,\gamma}$ and $\vartheta_{\sigma}=\vartheta_{\gamma}$, showing the result.
\hfill$\Box$

\begin{remark}
In the previous proposition, the condition that the $Q\mathbb{I}_{\Delta^c}$-invariant measure is finite is very important. Indeed, it is known that the statement of Proposition \ref{Prop-invariant} may be false even for a $\sigma$-finite non-null $Q\mathbb{I}_{\Delta^c}$-invariant measure: see Subsection 2.2.21 in \cite{piunovskiy13}.
\end{remark}

\section{Compactness of occupation measures}
\label{Sec-4}
We will first show in Theorem \ref{Eq-Characteristic++} that a measure is an occupation measure if and only if it satisfies the characteristic equation \eqref{eq-linear-equation} and an additional condition.
Roughly speaking, a measure satisfies this additional condition if it is absolutely continuous with respect to a reference probability measure~$\eta^\beta$ which will be introduced in Definition \ref{Def-lambda-beta}.
An interesting property that plays a central role in the analysis of the compactness of the set occupation measures is that the model $\bfrak{M}(\eta^\beta,\Delta)$ is also absorbing (see Proposition \ref{key-Proposition}).
Secondly, we will show in Theorem \ref{Relative-compactness-set-D} that the set of occupation measures is compact in the $ws$-topology if and only if the model is uniformly absorbing. 
In particular, this result highlights the fact that the assumption that the MDP is absorbing is not sufficient to guarantee the (relative) compactness of the set of occupation measures. 
This last point is illustrated by an example, at the end of this section, of a an absorbing model which is not uniformly absorbing, and whose  set of occupation measures $\bcal{O}_\eta$ is neither relatively compact nor compact.
As a corollary, since 
 a discounted MDP is uniformly absorbing, we recover the classical result stating that, under suitable continuity-compactness conditions, the set of occupation measures of a discounted model is compact.

\bigskip

\noindent\textbf{Condition (S)}
\begin{condit}
\item \label{Hyp-action} The action set $\AA$  is compact and the correspondence from $\XX$ to $\AA$ given by $x\mapsto \AA(x)$ is weakly measurable with nonempty compact values.
\item \label{Hyp-Q-strongly-continuous} For any $x\in \XX$ and $\Gamma\in \bfrak{X}$, the mapping $a\mapsto Q(\Gamma |x,\cdot)$ is continuous\footnote{In the previous version {\tt https://arxiv.org/abs/2309.07059v1} we assumed continuity of $a\mapsto Q(\Gamma |x,\cdot)$ on $\AA$.} on $\AA(x)$.
%For any $x\in \XX$ and $\Gamma\in \bfrak{X}$, the mapping $a\mapsto Q(\Gamma |x,\cdot)$ is continuous on $\AA$.
\end{condit}

\begin{lemma}
\label{Existence-xi-strongly-continuous}
Suppose that the Conditions \ref{Hyp-action}-\ref{Hyp-Q-strongly-continuous} are satisfied. There exists $\xi^*\in\mathbf{M}$ such that 
\begin{eqnarray}
Q(\cdot |x,a)\ll  Q_{\xi^*}(\cdot |x)\quad\hbox{for every $(x,a)\in\KK$}.
\label{Q-strongly-continuous}
\end{eqnarray}
\end{lemma}
\textbf{Proof.} The multifunction  from $\XX$ to~$\AA$ defined by $x\rightarrow \AA(x)$ is weakly measurable.
By \ref{Hyp-action} and Corollary 18.15 in \cite{aliprantis06}, we obtain the existence of a sequence $\{\xi_{n}\}_{n\in\NN}$ of measurable selectors for  the multifunction $x\mapsto \AA(x)$
satisfying 
\begin{equation}\label{eq-castaing}
\AA(x)= \widebar{\{\xi_{n}(x): n\in\NN\}}\quad\hbox{for any $x\in\XX$.}
\end{equation}
Define now  $\xi^*\in\mathbf{M}$ by means of 
$$\xi^*(da |x) = \sum_{k\in \NN} \frac{1}{2^{k+1}} \delta_{\xi_{k}(x)} (da)\quad\hbox{so that}\quad
Q_{\xi^*}(dy|x)=\sum_{k\in\NN} \frac{1}{2^{k+1}}Q(dy|x,\xi_k(x)).$$
To prove the result, fix arbitrary $(x,a)\in\KK$  and  $B\in\bfrak{X}$ such that $Q_{\xi^*}(B|x)=0$. This implies that $Q(B|x,\xi^*_k(x))=0$ for all $k\in\NN$.
Using \ref{Hyp-Q-strongly-continuous} and \eqref{eq-castaing} we obtain that $Q(B|x,a)=0$.
 \hfill$\Box$
 \\[10pt]\indent
 This result leads to the following definition
 \begin{definition}
\label{Def-lambda-beta}
The probability measure $\eta^\beta\in\bcal{P}(\XX)$ is
\begin{equation}\label{eq-def-lambda-beta}
\eta^\beta= (1-\beta) \sum_{k\in \NN} \beta^k \eta Q_{\xi^*}^k,
\end{equation}
where $\beta$ is some fixed parameter with $0<\beta<1$.
\end{definition}

Let us make some comments on this definition. By construction, the kernel $Q_{\xi^*}$ somehow approximates all possible one-step transitions (recall Lemma \ref{Existence-xi-strongly-continuous}). The measure $\eta^\beta$ spreads $Q_{\xi^*}$ over the  time horizon. Therefore, $\eta^\beta$ ``covers'' the whole MDP over space and time. That is why, as shown in our next result,  $\eta^\beta$ dominates any occupation measure.

\begin{lemma}
\label{Transition-kernel-absolute-continuity}
Suppose that $\bfrak{M}(\eta,\Delta)$ is absorbing and that the Conditions \ref{Hyp-action}-\ref{Hyp-Q-strongly-continuous}  are satisfied.
For every $\pi\in\mathbf{\Pi}$ we have
\begin{eqnarray*}
\mu^{\XX}_{\eta,\pi} \ll \mu^{\XX}_{\eta,\xi^*} \sim \eta^\beta_{\Delta^c} \ll \eta^\beta.
\end{eqnarray*}
\end{lemma}
\textbf{Proof:}
By Proposition \ref{prop-linear-equation-occupation}(ii), there exists $\sigma\in\mathbf{M}$ such that
\begin{equation}\label{eq-lambda-tool}
\mu_{\eta,\pi}^\XX=\eta \sum_{k=0}^{\infty} Q_{\sigma}^k\mathbb{I}_{\Delta^c}.
\end{equation}
We are going to show that for every $k\ge1$ and $x\in\XX$ we have $Q^k_\sigma(\cdot|x)\ll Q^k_{\xi^*}(\cdot|x)$. The proof is by induction. 

For the case $k=1$ suppose that $Q_{\xi^*}(B|x)=0$. By Lemma \ref{Existence-xi-strongly-continuous} this implies that $Q(B|x,a)=0$ for all $a\in\AA(x)$ and $Q_\sigma(B|x)=0$ follows. Assuming the result true for some $k$ and, for $k+1$, suppose that $Q^{k+1}_{\xi^*}(B|x)=0$.  Note that 
$$Q^{k+1}_{\xi^*}(B|x)=\int_\XX Q^k_{\xi^*}(B|y)Q_{\xi^*}(dy|x)=0$$
implies that $Q^k_{\xi^*}(B|y)=0$ for all $y\in C$ where $C$ is such that $Q_{\xi^*}(C|x)=1$. By the induction hypothesis, we have that 
$Q^k_\sigma(B|y)=0$ for all $y\in C$ with $Q_\sigma(C|x)=1$. This shows that, indeed, $Q^{k+1}_\sigma(B|x)=0$.

As a direct consequence we have that 
$$\eta\sum_{k=0}^\infty Q^k_\sigma \ll \eta\sum_{k=0}^\infty  Q^k_{\xi^*} \sim(1-\beta)\eta\sum_{k=0}^\infty \beta^k Q^k_{\xi^*}=\eta^\beta
$$
and the stated result follows from \eqref{eq-lambda-tool}.
\hfill$\Box$

\begin{proposition}
\label{key-Proposition}
Under Conditions \ref{Hyp-action}-\ref{Hyp-Q-strongly-continuous}, if the control model
$\bfrak{M}(\eta,\Delta)$ is absorbing then
$\bfrak{M}(\eta^\beta,\Delta)$ is also absorbing.
\end{proposition}
\textbf{Proof.} Consider an arbitrary $\pi\in\mathbf{\Pi}$. We have the following equalities
\begin{eqnarray*}
\mathbb{E}_{\eta^\beta,\pi}[T_\Delta] &=& \sum_{j=0}^\infty \mathbb{E}_{\eta^\beta,\pi}[\mathbf{I}_{\Delta^c}(X_j)]\\
&=& (1-\beta)\sum_{j=0}^\infty \sum_{k=0}^\infty \beta^k   \mathbb{E}_{\eta Q^k_{\xi^*},\pi}[\mathbf{I}_{\Delta^c}(X_j)].
\end{eqnarray*}
Observe that using the policy $\pi$ for the initial distribution $\eta Q^k_{\xi^*}$ is equivalent to using the policy $\gamma^k\in\mathbf{\Pi}$ given by
$$\gamma^k_j(da|x_0,a_0,\ldots,x_j)=\xi^*(da|x_j)\quad\hbox{for $0\le j<k$}$$
and 
$$\gamma^k_j(da|x_0,a_0,\ldots,x_j)=\pi_{j-k}(da|x_{k},a_{k},\ldots,x_{j})\quad\hbox{for $j\ge k$}$$
for the initial distribution $\eta$ just by making a  shift of $k$ time units.
Therefore, 
$$ \mathbb{E}_{\eta Q^k_{\xi^*},\pi}[\mathbf{I}_{\Delta^c}(X_j)]= \mathbb{E}_{\eta ,\gamma^k}[\mathbf{I}_{\Delta^c}(X_{j+k})].
$$
It follows that 
\begin{eqnarray*}
\mathbb{E}_{\eta^\beta,\pi}[T_\Delta] &=& (1-\beta)\sum_{j=0}^\infty \sum_{k=0}^\infty \beta^k  
 \mathbb{E}_{\eta ,\gamma^k}[\mathbf{I}_{\Delta^c}(X_{j+k})]\\
&=& 
 (1-\beta)\sum_{k=0}^\infty \beta^k   \sum_{j=0}^\infty 
 \mathbb{E}_{\eta ,\gamma^k}[\mathbf{I}_{\Delta^c}(X_{j+k})]\\
 &\le& (1-\beta)\sum_{k=0}^\infty \beta^k  \mathbb{E}_{\eta ,\gamma^k}[T_\Delta].
 \end{eqnarray*}
 The model $\bfrak{M}(\eta,\Delta)$ being absorbing, we have that 
 $$\sup_{\pi\in\mathbf{\Pi}} \mathbb{E}_{\eta ,\pi}[T_\Delta]=\mathbf{c}<\infty$$
 and so
 $\mathbb{E}_{\eta^\beta,\pi}[T_\Delta] \le \mathbf{c}$, which shows that $\bfrak{M}(\eta^\beta,\Delta)$ is absorbing as well.
 \hfill$\Box$
 %\\[10pt]

\begin{proposition}
\label{prop-survival-TDelta}
Suppose  that $\bfrak{M}(\eta,\Delta)$ is absorbing and that the Conditions \ref{Hyp-action}-\ref{Hyp-Q-strongly-continuous}  are satisfied.
Let $\boldsymbol{\Gamma}$ be an arbitrary subset of $\mathbf{M}$ and let $\{h_{\pi}\}_{\pi\in\boldsymbol{\Gamma}}$ be a family of non-negative functions in $L^{1}(\XX,\bfrak{X},\eta^\beta)$ which are uniformly $\eta^\beta$-integrable. Under these conditions,
$$
\lim_{t\rightarrow\infty} \sup_{\pi\in \boldsymbol{\Gamma}} \int_\XX Q^t_\pi(\Delta^c|x)h_\pi(x)\eta^\beta(dx)= 0.
$$
\end{proposition}
\textbf{Proof.}
Consider a fixed arbitrary $\epsilon>0$. By the uniform integrability hypothesis, there exists $c_{\epsilon}>0$ such that
$$ \sup_{\pi\in \boldsymbol{\Gamma}} \int_{\{x\in\XX: h_{\pi}(x)>c_{\epsilon}\}} h_{\pi}(x)  \eta^\beta(dx)\leq \epsilon.$$
Therefore, for any  $\pi\in\mathbf{\Gamma}$  and $t\ge1$
$$
 \int_\XX Q^t_\pi(\Delta^c|x)h_\pi(x) \eta^\beta(dx)
 \leq \epsilon+c_{\epsilon} \mathbb{P}_{\eta^\beta,\pi}\{T_{\Delta}>t\}
\leq \epsilon+ \frac{c_{\epsilon}}{t} \cdot\mathbb{E}_{\eta^\beta,\pi}[T_{\Delta}]
\leq \epsilon+ \frac{c_{\epsilon}}{t} \sup_{\pi\in \mathbf{\Pi}} \mathbb{E}_{\eta^\beta,\pi}[T_{\Delta}].
$$
By Proposition \ref{key-Proposition} we have that the above supremum is finite.
Hence, for~$t$ sufficiently large we obtain that 
$$\sup_{\pi\in \boldsymbol{\Gamma}} \int_\XX Q^t_\pi(\Delta^c|x)h_\pi(x)\eta^\beta(dx)<2\epsilon$$
and the result follows.
$\hfill$ $\Box$\\[10pt]
\indent
Our next result gives a characterization of $\bcal{O}_\eta$ based on the probability measure $\eta^\beta$.
\begin{theorem}\label{Eq-Characteristic++}
Let $\bfrak{M}(\eta,\Delta)$ be an absorbing model satisfying the conditions  \ref{Hyp-action}-\ref{Hyp-Q-strongly-continuous}. 
Given a measure  $\mu\in\bcal{C}_\eta$, we have $\mu\in\bcal{O}_\eta$ if and only if $\mu^{\XX} \ll\eta^\beta$. Equivalently,
$$\mathcal{O_{\eta}}=\Big\{ \mu\in \MM_{+}(\XX\times\AA) : \mu(\KK^c)=0; \quad  \mu^\XX=\big(\eta+\mu^\XX Q_{\sigma}\big) \mathbb{I}_{\Delta^c}\quad\hbox{and}\quad \mu^{\XX} \ll\eta^\beta \Big\}.$$
\end{theorem}
\textbf{Proof.} Clearly, if $\mu\in \bcal{O}_\eta$ then $\mu^\XX \ll\eta^\beta$ as a direct consequence of Lemma \ref{Transition-kernel-absolute-continuity}. 

Conversely, let us assume that $\mu\in\bcal{C}_\eta$ is such that $\mu^\XX\ll\eta^\beta$. 
Using Theorem~\ref{prop-characterizations}(i), we can find $\sigma\in\mathbf{M}$ and a $Q\mathbb{I}_{\Delta^c}$-invariant measure $\vartheta\in\MM^+(\XX\times\AA)$ for the kernel $Q\mathbb{I}_{\Delta^c}$ such that   
$$
\mu=\mu_{\eta,\sigma}+\vartheta
\quad\hbox{and}\quad
\vartheta^\XX=\vartheta^\XX Q_{\sigma}\mathbb{I}_{\Delta^c}.$$
This implies that
$\vartheta^\XX\ll\eta^\beta$ and so for every $t\ge1$
\begin{align*}
\vartheta^\XX(\XX)=\vartheta^\XX (Q_{\sigma}\mathbb{I}_{\Delta^c})^t(\XX)=\vartheta^\XX Q^t_{\sigma}(\Delta^c)=\int_{\XX} Q^t_{\sigma}(\Delta^c |x) \frac{d\vartheta^\XX}{d\eta^\beta}(x) \eta^\beta(dx).
\end{align*}
Applying Proposition \ref{prop-survival-TDelta} to the set $\mathbf{\Gamma}=\{\sigma\}$ and the function $d\vartheta^\XX/d\eta^\beta$, we can take the limit as $t\rightarrow\infty$ in the last expression to obtain that   $\vartheta^\XX(\XX)=0$. This shows that, indeed, $\mu\in\bcal{O}_\eta$.
\hfill $\Box$

\begin{lemma}
\label{Balder-lsc}
Suppose $\bfrak{M}(\eta,\Delta)$ is absorbing and that the Conditions \ref{Hyp-action}-\ref{Hyp-Q-strongly-continuous} hold.
Let $\Psi$ be a measurable non-negative function defined on $\XX\times\AA$ such that $\Psi(x,\cdot)$ is lower semi-continuous on~$\AA$ for each $x\in\XX$. 
Then there exist $\mathcal{N}\in\bfrak{X}$ and a sequence $\{\Psi_n\}_{n\in\NN}$ of  $\bfrak{X}\otimes\bfrak{B}(\mathbf{A})$-measurable bounded Carath\'eodory functions satisfying
$\Psi_n\uparrow\mathbf{I}_{(\mathcal{N}^c\times\AA)} \Psi$
with $\eta^\beta(\mathcal{N})=0$.
\end{lemma}
\textbf{Proof.} Let us denote by $\bar{\bfrak{X}}$ the completion of $\bfrak{X}$ with respect to $\lambda_{\beta}$.
By \mbox{\cite[Proposition 3.9]{florescu12}},
 there exists a sequence $\{\bar{\Psi}_n\}_{n\in\NN}$ of  $\bar{\bfrak{X}}\otimes\bfrak{B}(\mathbf{A})$-measurable bounded Carath\'eodory functions
such that ${\bar{\Psi}}_{n} \uparrow  \Psi$ as $n\rightarrow\infty$. 
For every $n\in\NN$, applying \cite[Lemma A-6]{balder84}  to $\bar{\Psi}_n$ (being both a l.s.c. and a u.s.c. integrand), there exists a set $N_n\in\bfrak{X}$ with $\eta^\beta(N_n)=0$ and a $\bfrak{X}\otimes\bfrak{B}(\mathbf{A})$-measurable bounded Carath\'eodory function $\tilde{\Psi}_{n}$ such that 
$$\tilde{\Psi}_n(x,a)=\bar{\Psi}_n(x,a)\quad\hbox{for all $x\in\XX-N_n$ and $a\in\AA$}.$$
Letting $\mathcal{N}=\cup_{n\in\NN} N_n$, we have $\eta^\beta(\mathcal{N})=0$. Let us define $\Psi_n(x,a)=\mathbf{I}_{\mathcal{N}^c}(x)\cdot \tilde\Psi_n(x,a)$ for
$(x,a)$ in $\XX\times\AA$. It should be clear that each $\Psi_n$ is a $\bfrak{X}\otimes\bfrak{B}(\mathbf{A})$-measurable bounded Carath\'eodory function satisfying the required property.
\hfill$\Box$

\begin{lemma}
\label{absolute-cont+Kc} 
Suppose $\bfrak{M}(\eta,\Delta)$ is absorbing and that the Conditions \ref{Hyp-action}-\ref{Hyp-Q-strongly-continuous} hold.
Let us consider a net $\{\nu_{\alpha}\}_{\alpha\in\mathcal{I}}$ in $\bcal{M}^+(\XX\times\AA)$ converging to some $\nu\in\bcal{M}^+(\XX\times\AA)$ in the $ws$-topology.
Assume that $\nu^{\XX}_{\alpha} \ll\eta^\beta$ and $\nu_{\alpha}(\KK^c)=0$ for any $\alpha\in\mathcal{I}$.
Then $\nu^{\XX} \ll\eta^\beta$ and $\nu(\KK^c)=0$
\end{lemma}
\textbf{Proof.} 
Let $B\in\bfrak{X}$ such that $\eta^\beta(B)=0$. Then by hypothesis, $\nu^{\XX}_{\alpha}(B)=0$ for any $\alpha\in\mathcal{I}$ and so,
$\nu^{\XX}(B)=\lim_{\alpha}\nu^{\XX}_{\alpha}(B)=0$ giving the first part of the result.
Let us show that we also have $\nu(\KK^c)=0$.
Since the action sets $\AA(x)$ are compact, we have that the function $(x,a)\mapsto \mathbf{I}_{\KK^c}(x,a)$ is measurable, non-negative and lower semi-continuous on $\AA$ for each $x\in\XX$. 
From Lemma \ref{Balder-lsc}, there exist $\mathcal{N}\in\bfrak{X}$ and a sequence $\{\Psi_n\}_{n\in\NN}$ of  $\bfrak{X}\otimes\bfrak{B}(\mathbf{A})$-measurable bounded Carath\'eodory functions satisfying
$\Psi_n\uparrow\mathbf{I}_{(\mathcal{N}^c\times\AA)\cap\KK^c}$ with $\eta^\beta(\mathcal{N})=0$.
Hence, given $\epsilon>0$ we can find some $n_0$ with $\ds \nu\big(\KK^c\big)\le\int_{\XX\times\AA} \Psi_{n_0}d\nu+\epsilon/2$
since $\nu^\XX(\mathcal{N})=0$. 
There exists some $\alpha_0$ such that for every $\alpha\succ\alpha_0$
$$\nu(\KK^c)\le\int_{\XX\times\AA} \Psi_{n_0}d\nu+\epsilon/2\le \int_{\XX\times\AA} \Psi_{n_0}d\nu_\alpha+\epsilon
\le \nu_\alpha(\KK^c\big)+\epsilon=\epsilon.$$
We conclude that $\nu(\KK^c)=0$. 
\hfill$\Box$
\\[10pt]\indent
The next lemma shows, roughly speaking, that for a measurable set $B$ in $\XX$, the restriction of the function $Q(B|\cdot)$ on $\KK$ can be extended on $\mathcal{N}^c\times\AA$
to a Carath\'eodory function defined $\XX\times\AA$ for $\mathcal{N}\in\bfrak{X}$ such that $\eta^\beta(\mathcal{N})=0$.
\begin{lemma}
Suppose that the Conditions \ref{Hyp-action}-\ref{Hyp-Q-strongly-continuous} are satisfied. 
\label{caratheodory}
For each measurable set $B\in\bfrak{X}$ there exist $\mathcal{N}\in\bfrak{X}$ and a Carath\'eodory function $\widehat{Q}(B| \cdot,\cdot)$ on $\XX\times\AA$ satisfying
\begin{align*}
\widehat{Q}(B|x,a) = Q(B|x,a)
\end{align*}
for every $x\in\mathcal{N}^c$ and $a\in\AA$, with $\eta^\beta(\mathcal{N})=0$.
\end{lemma}
\textbf{Proof.}
Let us write $Q\big|_{\KK}(B| \cdot,\cdot)$ for restriction of $Q(B| \cdot,\cdot)$ to $\KK$ and
denote by $\widebar{\bfrak{X}}$ the completion of $\bfrak{X}$ with respect to $\lambda_{\beta}$.
Let $B\in\bfrak{X}$ be a fixed set.
By using Corollary 3 in \cite{kucia98} there exists a $\widebar{\bfrak{X}}\otimes\bfrak{B}(\mathbf{A})$-measurable Carath\'eodory extension of $Q\big|_{\KK}(B|\cdot,\cdot)$ on $\XX\times\AA$
denoted by $\widebar{Q}(B|\cdot,\cdot)$.
Applying Lemma A-6 in \cite{balder84}  to $\widebar{Q}(B|\cdot,\cdot)$ (being both a l.s.c. and a u.s.c. integrand), there exists a measurable set $\mathcal{N}\in\bfrak{X}$ with
$\eta^\beta(\mathcal{N})=0$ and a $\bfrak{X}\otimes\bfrak{B}(\mathbf{A})$-measurable Carath\'eodory function $\widetilde{Q}(B|\cdot,\cdot)$ such that 
$$\widebar{Q}(B|x,a)=\widetilde{Q}(B|x,a)\quad\hbox{for all $x\in\mathcal{N}^c$ and $a\in\AA$}.$$
Finally, the function $\widehat{Q}(B|\cdot,\cdot)$ defined on $\XX\times\AA$ by $\widehat{Q}(B|x,a)=\mathbf{I}_{\mathcal{N}^c}(x)\widetilde{Q}(B|x,a)$ satisfies the desired properties.
 \hfill$\Box$

\begin{lemma}
\label{Relatively-s-compactness-Q}
Assume that the Conditions \ref{Hyp-action}-\ref{Hyp-Q-strongly-continuous} hold.
Let $\{\nu_{\pi}\}_{\pi\in\mathbf{M}}$ be a relatively $s$-compact subset of $\bcal{M}^+(\XX)$ such that $\nu_{\pi} \ll\eta^\beta$ for any $\pi\in\mathbf{M}$.
Then set $\{\nu_{\pi} Q_{\pi}\}_{\pi\in\mathbf{M}}$ is also relatively $s$-compact subset of
$\bcal{M}^+(\XX)$.
\end{lemma}
\textbf{Proof.} 
Define $$\mathbf{\Lambda}=\{\nu_{\pi}\otimes \pi\}_{\pi\in\mathbf{M}}\subseteq\bcal{M}^+(\XX\times\AA).$$
By hypothesis,  $\mathbf{\Lambda}^\XX$ is relatively $s$-compact and $\mathbf{\Lambda}^\AA$ is tight since $\AA$ is a compact metric space. 
Therefore, $\mathbf{\Lambda}$ is relatively sequentially $ws$-compact by using Theorem 2.5 in \cite{balder01}. So, there exists a sequence $\{\pi_n\}$ in $\mathbf{M}$ such that 
\begin{equation}\label{eq-relatively-ws}
\{\nu_{\pi_{n}}\otimes \pi_{n}\}\stackrel{ws}{\longrightarrow}  \mu
\end{equation}
 for some measure $\mu\in\bcal{M}^+(\XX\times\AA)$.
Observe that  $\nu_{\pi_{n}}\otimes \pi_{n}(\KK^{c})=0$ and by Lemma \ref{absolute-cont+Kc}, it follows also that  $\mu(\KK^{c})=0$.
For $B\in\bfrak{X}$, let us consider the function $\widehat{Q}(B| \cdot,\cdot)$ defined on $\XX\times\AA$ as introduced in Lemma \ref{caratheodory}.
From equation \eqref{eq-relatively-ws}, we obtain
$$\nu_{\pi_n} Q_{\pi_n}(B)=\int_{\XX\times\AA} \widehat{Q}(B|x,a)d(\nu_{\pi_{n}}\otimes \pi_{n})
\rightarrow \int_{\XX\times\AA}\widehat{Q}(B|x,a)d\mu= \mu Q(B),$$
which shows that 
$$\nu_{\pi_n}Q_{\pi_n}\stackrel{s}{\rightarrow}\mu Q,$$
establishing that $\{\nu_{\pi} Q_{\pi}\}_{\pi\in\mathbf{M}}$ is relatively sequentially $s$-compact.
Finally, \cite[Proposition 2.2]{balder01} or \cite[Theorem 2.6]{ganssler71} give the result on relative $s$-compactness.
\hfill$\Box$
\\[10pt]\indent
Given the initial distribution $\eta\in\PP(\XX)$, a  policy $\pi\in\mathbf{\Pi}$, and $k\in\NN$ we define $\mu_{\eta,\pi,k}$ as the distribution of $(X_k,A_k)$ under  $\mathbb{P}_{\eta,\pi}$. In symbols, we let
$$\mu_{\eta,\pi,k}=\mathbb{P}_{\eta,\pi}\circ (X_k,A_k)^{-1}\in\PP(\XX\times\AA).$$
We consider the following subsets of measures
$$\bcal{D}_{\eta,k}^s=\{\mu_{\eta,\pi,k}\}_{\pi\in\mathbf{M}}\quad\hbox{for $k\in\NN$},$$
which correspond to the families of state-action distributions of the process at time $k$ for stationary policies.

\begin{lemma}
\label{Relatively-ws-compactness-Ds}
Assume that the Conditions \ref{Hyp-action}-\ref{Hyp-Q-strongly-continuous} hold.
For each $k\in\NN$, the sets $ \bcal{D}_{\eta,k}^s$ are relatively $ws$-compact and 
\begin{equation}
\label{eq-to-show-relatively-compact}
\lim_{n\rightarrow\infty} \sup_{\pi\in\mathbf{M}} \mu^\XX_{\eta,\pi,k}(\Gamma_{n})=0
\end{equation}
for any decreasing sequence of sets $\Gamma_{n}\in\bfrak{X}$ with $\Gamma_{n}\downarrow\emptyset$.
\end{lemma}
\textbf{Proof.} First of all we prove by induction that $ (\bcal{D}_{\eta,k}^s)^\XX$ is relatively $s$-compact in $\MM^+(\XX)$ for every $k\in\NN$. For $k=0$ we have $  (\bcal{D}_{\eta,k}^s)^\XX=\{\eta\}$ and the result is obvious.
Assuming the result is true for some $k\ge0$, observe that 
$\mu^\XX_{\eta,\pi,k+1}=\mu^\XX_{\eta,\pi,k}Q_{\pi}$ and $\mu^\XX_{\eta,\pi,k}\ll\eta^{\beta}$ for any $\pi\in\mathbf{M}$.
We  get the claim by applying Lemma \ref{Relatively-s-compactness-Q}. 

Once we know that $\{\mu_{\eta,\pi,k}^\XX\}_{\pi\in\mathbf{M}}$ is relatively $s$-compact and 
$\AA$ being a  compact metric space, it follows that
$$\{\mu^\XX_{\eta,\pi,k}\otimes\pi\}_{\pi\in\mathbf{M}}=\{\mu_{\eta,\pi,k}\}_{\pi\in\mathbf{M}}=\bcal{D}^s_{\eta,k}$$
 is relatively $ws$-compact by Theorem 5.2 in \cite{balder01}, showing the first claim. The last part of the result is obtained
 by using Lemma 4.6.5 and Theorem 4.7.25 in~\cite{bogachev07}.
 \hfill $\Box$
 
\begin{lemma}
\label{Closedness-set-Dm} 
Suppose $\bfrak{M}(\eta,\Delta)$ is absorbing and that the Conditions \ref{Hyp-action}-\ref{Hyp-Q-strongly-continuous} hold.
The set $\bcal{O}_{\eta}$ is  $ws$-closed in $\bcal{M}^{+}(\XX\times\AA)$.
\end{lemma}
\textbf{Proof.} 
Let us consider a net $\{\nu_{\alpha}\}_{\alpha\in\mathcal{I}}$ in $\bcal{O}_{\eta}$ converging to some $\nu\in\bcal{M}^+(\XX\times\AA)$ in the $ws$-topology.
From Theorem~\ref{Eq-Characteristic++} we have that for any $\alpha\in\mathcal{I}$
\begin{align*}
\nu_{\alpha}(\KK^c)=0, \qquad \nu^\XX_{\alpha}=(\eta+  \nu_{\alpha} Q )\mathbb{I}_{\Delta^c},
\quad\hbox{and}\quad
\nu^{\XX}_{\alpha} \ll\eta^\beta.
\end{align*}
It easily follows by Lemma \ref{absolute-cont+Kc} that that the limiting measure $\nu$ satisfies $\nu^\XX\ll\eta^\beta$  and  $\nu(\KK^c)=0$.
Let us show that we also have $\nu^\XX=(\eta+  \nu Q )\mathbb{I}_{\Delta^c}$.
For $B\in\bfrak{X}$, let us consider the function $\widehat{Q}(B| \cdot,\cdot)$ defined on $\XX\times\AA$ as introduced in Lemma \ref{caratheodory}.
Therefore,
\begin{align*}
\lim_{\alpha}  \nu_{\alpha} Q \mathbb{I}_{\Delta^c}(B) & =  \lim_{\alpha} \int_{\XX\times\AA} \mathbf{I}_{\Delta^c}(x) \widehat{Q}(B|x,a) \nu_{\alpha}(dx,da)
=  \int_{\XX\times\AA} \mathbf{I}_{\Delta^c}(x) \widehat{Q}(B|x,a) \nu(dx,da) \\
& = \nu Q \mathbb{I}_{\Delta^c}(B)
\end{align*}
showing that $\nu^\XX=(\eta+  \nu Q )\mathbb{I}_{\Delta^c}$.
This implies that $\nu\in\bcal{O}_\eta$ and, therefore, $\bcal{O}_\eta$ is closed.
\hfill$\Box$

\begin{theorem}
\label{Relative-compactness-set-D} 
Suppose $\bfrak{M}(\eta,\Delta)$ is absorbing and that the Conditions \ref{Hyp-action}-\ref{Hyp-Q-strongly-continuous} are satisfied. The following statements are equivalent.
\begin{itemize}
\item[$(a)$] $\bfrak{M}(\eta,\Delta)$ is uniformly absorbing.
\item[$(b)$] The set $\bcal{O}_\eta$ is relatively $ws$-compact.
\item[$(c)$] The set $\bcal{O}_\eta$ is $ws$-compact.
\end{itemize}
\end{theorem}
\textbf{Proof.} 
$(a)\Rightarrow (b)$.
The set of $\AA$-marginal measures of $\bcal{O}_\eta$ is tight since $\AA$ is compact.
Recalling that $\bcal{O}_\eta=\bcal{O}_\eta^s$, by Proposition \ref{prop-linear-equation-occupation}(ii) and
applying  \cite[ Theorem 5.2]{balder01} we need to show that the set of $\XX$-marginal measures of $\bcal{O}^s_\eta$ is relatively $s$-compact or, equivalently, that 
\begin{equation}\label{eq-to-show-relatively-compact-bis}
\lim_{n\rightarrow\infty} \sup_{\pi\in\mathbf{M}} \mu^\XX_{\eta,\pi}(\Gamma_{n})=0
\end{equation}
for any decreasing sequence of sets $\Gamma_{n}\in\bfrak{X}$ with $\Gamma_{n}\downarrow\emptyset$
(see Lemma 4.6.5 and Theorem 4.7.25~in~\cite{bogachev07}).
Now, observe that for every $k\in\NN$
\begin{eqnarray}
\sup_{\pi\in\mathbf{M}} \mu^\XX_{\eta,\pi}(\Gamma_n)  &\le&  \sum_{t=0}^k  \sup_{\pi\in\mathbf{M}} \mu^\XX_{\eta,\pi,t}(\Gamma_{n})
+\sup_{\pi\in\mathbf{M}} \sum_{t>k}\mathbb{P}_{\eta,\pi_\mu}\{T_\Delta>t\}.
\label{eq-tool-relatively-compact}
\end{eqnarray}
We have that 
$\sup_{\pi\in\mathbf{M}} \sum_{t>k}\mathbb{P}_{\eta,\pi_\mu}\{T_\Delta>t\}$ as $k\rightarrow\infty$ since 
$\bfrak{M}(\eta,\Delta)$ is uniformly absorbing by hypothesis.
%Moreover, the leftmost term in \eqref{eq-tool-relatively-compact} converges to zero as $k\rightarrow\infty$
Morerover, for each $k \in\mathbb{N}$, the leftmost term on the right-hand side of \eqref{eq-tool-relatively-compact} converges to zero as $n\rightarrow\infty$
by using Lemma \ref{Relatively-ws-compactness-Ds} showing that the limit \eqref{eq-to-show-relatively-compact-bis} holds.
This shows that $\bcal{O}_\eta$ is relatively compact for the $ws$-topology. 
\\[5pt]
\noindent
$(b)\Rightarrow (a)$.
Since $\bcal{O}_\eta$ is relatively compact for the $ws$-topology, it follows from  \cite[Theorem 5.2]{balder01} that the set of $\XX$-marginal measures of $\bcal{O}_\eta=\bcal{O}_\eta^s$ is relatively s-compact.
Recalling that  $\mu^\XX_{\eta,\pi}\ll \eta^\beta$  for any $\pi\in\mathbf{M}$, using  Lemma \ref{Transition-kernel-absolute-continuity},
Proposition 2.2 in \cite{balder01}, and Corollary 2.7 in \cite{ganssler71}, we get that the family $\{h_\pi\}_{\pi\in\mathbf{M}}$ of density functions 
$h_\pi=d\mu^\XX_{\eta,\pi}/d\eta^\beta$ is uniformly $\eta^\beta$-integrable.
Now, observe that for $\pi\in\mathbf{M}$,
$$ \sum_{k=t}^{\infty} \mathbb{P}_{\eta,\pi}\{T_\Delta>k\}=
\sum_{j=0}^\infty \mathbb{P}_{\eta Q_\pi^j,\pi}\{T_\Delta>t\}
=
 \mu^\XX_{\eta,\pi} Q_{\pi}^{t}(\Delta^c)=\int_{\XX} Q_{\pi}^{t}(\Delta^c |x) h_\pi(x) \eta^\beta(dx)$$
and by using Proposition \ref{prop-survival-TDelta} we can conclude that the rightmost term in the previous equation converges to zero uniformly in $\pi\in\mathbf{M}$ as $t\rightarrow\infty$.
This establishes that $\bfrak{M}(\eta,\Delta)$ is indeed uniformly absorbing.
\\[5pt]
\indent
$(b)\Leftrightarrow (c)$. This is obvious using Lemma \ref{Closedness-set-Dm} and recalling that $\bcal{O}_\eta=\bcal{O}_\eta^s$ 
\hfill$\Box$

\begin{example}
In this example, we show that for a specific model which is not uniformly absorbing, the set of occupation measures $\bcal{O}_\eta$ is not relatively compact and so, not compact.
Consider the model defined with state and actions spaces
$$\mathbf{X}=\{ (i,j)\in\NN\times\NN: 0\le j\le 2^i-1\} \cup\{\bar{x}\}\quad\hbox{and}\quad\AA=\{c,s\},$$
with 
$\AA(x)=\{c,s\}$ for $x\in\NN\times\{0\}$ and $\AA(x)=\{s\}$ otherwise.
The transition kernel is given by:
\begin{itemize}
\item $Q\big(\bar{x}\mid (i,0),c)\big)=\frac{1}{2}$ and $Q\big((i+1,0)\mid (i,0),c)\big)=\frac{1}{2}$ for  $i\in\NN$;
\item $Q\big((i,j+1)\mid (i,j),s\big)=1$ for $i\in\NN$ and $0\le j<2^i-1$;
\item $Q\big(\bar{x}\mid (i,2^i-1),s\big)=1$ for $i\in\NN$;
\item $Q(\bar{x}\mid \bar{x},s)=1$.
\end{itemize}
In \cite[Example 3.13]{piunovskiy19}, it has been shown that this model is absorbing to $\Delta=\{\bar{x}\}$ with initial distribution
$\eta=\delta_{(0,0)}$ but not uniformly absorbing.
Our objective is to illustrate by means of this example that the set of occupation measures $\bcal{O}_{\eta}$ is not compact for the $ws$-topology, as can be derived from Theorem \ref{Relative-compactness-set-D}.

Let us consider the sequence $\{\gamma_{t}\}_{t\in\NN}$ of nonrandomized stationary policies in $\mathbf{M}$ defined as follows. The policy $\gamma_t(da|x)$ takes the action $c$ for all states $(i,0)$ for $0\le i<t$ and takes the action~$s$ for all states $(i,0)$ with $i\ge t$. 
If $\bcal{O}_\eta$  were relatively  compact, then  for any decreasing sequence of sets  $\Gamma_{n}\downarrow\emptyset$ we would have
\begin{equation}\label{example-contradiction}
\lim_{n\rightarrow\infty} \sup_{t\in\NN} \mu^\XX_{\eta,\gamma_{t}}(\Gamma_{n})=0.
\end{equation}
Let us consider
$$\Gamma_{n}=\big\{(i,j)\in\XX: i\geq n \text{ and }1\le  j\le 2^i-1\}\big\}$$
which indeed satisfies  $\Gamma_{n}\downarrow\emptyset$. Observe that $$\mathbb{P}_{\eta,\gamma_{n}} \{ X_{n+j} = (n,j)\}=\frac{1}{2^n}\quad\hbox {for $n\geq 1$
and $j\in\{0,\cdots,2^n-1\}$}.$$
Therefore,
$$ \sup_{t\in\NN} \mu^\XX_{\eta,\gamma_{t}}(\Gamma_{n}) \geq \mu^\XX_{\eta,\gamma_{n}}(\Gamma_{n}) =2^n \cdot\frac{1}{2^n} =1,$$
contradicting \eqref{example-contradiction} implying that  $\bcal{O}_\eta$  is not relatively compact
This exhibits that, for this control model which is not uniformly absorbing, the set of occupation measures $\bcal{O}_\eta$ is not compact.    
\end{example}

Note that, in this example, the weak topology and the $ws$-topology on $\bcal{P}(\XX\times\AA)$ are the same. In view of this example, the result stated in \cite[Lemma 4.7]{feinberg12} is not accurate since, in order to have compactness of occupation measures, the absorbing property does not suffice, while the  uniformly absorbing condition is precisely the necessary and sufficient condition.


\begin{thebibliography}{10}

\bibitem{aliprantis06}
C.D. Aliprantis and K.C. Border.
\newblock {\em Infinite dimensional analysis}.
\newblock Springer, Berlin, 2006.

\bibitem{altman99}
E.~Altman.
\newblock {\em Constrained {M}arkov decision processes}.
\newblock Stochastic Modeling. Chapman \& Hall/CRC, Boca Raton, FL, 1999.

\bibitem{balder84}
Erik~J. Balder.
\newblock A general approach to lower semicontinuity and lower closure in
  optimal control theory.
\newblock {\em SIAM J. Control Optim.}, 22(4):570--598, 1984.

\bibitem{balder01}
Erik~J. Balder.
\newblock On ws-convergence of product measures.
\newblock {\em Math. Oper. Res.}, 26(3):494--518, 2001.

\bibitem{bogachev07}
V.I. Bogachev.
\newblock {\em Measure theory. {V}ol. {I}, {II}}.
\newblock Springer-Verlag, Berlin, 2007.

\bibitem{borkar02}
Vivek~S. Borkar.
\newblock Convex analytic methods in {M}arkov decision processes.
\newblock In {\em Handbook of {M}arkov decision processes}, volume~40 of {\em
  Internat. Ser. Oper. Res. Management Sci.}, pages 347--375. Kluwer Acad.
  Publ., Boston, MA, 2002.

\bibitem{dufour23}
Fran\c{c}ois Dufour and Tom\'{a}s Prieto-Rumeau.
\newblock Nash equilibria for total expected reward absorbing {M}arkov games: the constrained and unconstrained cases.
{\em Appl. Math. Optim.}, 89(2): paper No. 34, 2024.

\bibitem{dynkin79}
E.B. Dynkin and A.A. Yushkevich.
\newblock {\em Controlled {M}arkov processes}, volume 235 of {\em Grundlehren
  der Mathematischen Wissenschaften}.
\newblock Springer-Verlag, Berlin, 1979.

\bibitem{piunovskiy19}
E.A. Feinberg and A.~Piunovskiy.
\newblock Sufficiency of deterministic policies for atomless discounted and
  uniformly absorbing {MDP}s with multiple criteria.
\newblock {\em SIAM J. Control Optim.}, 57(1):163--191, 2019.

\bibitem{feinberg12}
E.A. Feinberg and U.G. Rothblum.
\newblock Splitting randomized stationary policies in total-reward {M}arkov
  decision processes.
\newblock {\em Math. Oper. Res.}, 37(1):129--153, 2012.

\bibitem{florescu12}
L.C. Florescu and C.~Godet-Thobie.
\newblock {\em Young measures and compactness in measure spaces}.
\newblock De Gruyter, Berlin, 2012.

\bibitem{ganssler71}
Peter G\"{a}nssler.
\newblock Compactness and sequential compactness in spaces of measures.
\newblock {\em Z. Wahrscheinlichkeitstheorie und Verw. Gebiete}, 17:124--146,
  1971.

\bibitem{he-sun17}
Wei He and Yeneng Sun.
\newblock Stationary {M}arkov perfect equilibria in discounted stochastic
  games.
\newblock {\em J. Econom. Theory}, 169:35--61, 2017.

\bibitem{hernandez96}
On\'esimo Hern{\'a}ndez-Lerma and Jean-Bernard Lasserre.
\newblock {\em Discrete-time {M}arkov control processes: Basic optimality
  criteria}, volume~30 of {\em Applications of Mathematics}.
\newblock Springer-Verlag, New York, 1996.

\bibitem{hernandez99}
On\'esimo Hern{\'a}ndez-Lerma and Jean-Bernard Lasserre.
\newblock {\em Further topics on discrete-time {M}arkov control processes},
  volume~42 of {\em Applications of Mathematics}.
\newblock Springer-Verlag, New York, 1999.

\bibitem{kucia98}
A.~Kucia.
\newblock Some results on {C}arath\'{e}odory selections and extensions.
\newblock {\em J. Math. Anal. Appl.}, 223(1):302--318, 1998.

\bibitem{neveu70}
Jacques Neveu.
\newblock {\em Bases math\'{e}matiques du calcul des probabilit\'{e}s}.
\newblock Masson et Cie, \'{E}diteurs, Paris, 1970.
\newblock Pr\'{e}face de R. Fortet, Deuxi\`eme \'{e}dition, revue et
  corrig\'{e}e.

\bibitem{nowak92}
Andrzej~S. Nowak and Tirukkannamangai E.~S. Raghavan.
\newblock Existence of stationary correlated equilibria with symmetric
  information for discounted stochastic games.
\newblock {\em Math. Oper. Res.}, 17(3):519--526, 1992.

\bibitem{piunovskiy97}
A.~B. Piunovskiy.
\newblock {\em Optimal control of random sequences in problems with
  constraints}, volume 410 of {\em Mathematics and its Applications}.
\newblock Kluwer Academic Publishers, Dordrecht, 1997.
\newblock With a preface by V. B. Kolmanovskii and A. N. Shiryaev.

\bibitem{piunovskiy13}
A.~B. Piunovskiy.
\newblock {\em Examples in {M}arkov decision processes}, volume~2 of {\em
  Imperial College Press Optimization Series}.
\newblock Imperial College Press, London, 2013.

\bibitem{schal75}
Manfred Sch{\"a}l.
\newblock On dynamic programming: compactness of the space of policies.
\newblock {\em Stochastic Processes Appl.}, 3(4):345--364, 1975.

\bibitem{valadier73}
Michel Valadier.
\newblock D\'{e}sint\'{e}gration d'une mesure sur un produit.
\newblock {\em C. R. Acad. Sci. Paris S\'{e}r. A}, 276:33--35, 1973.

\end{thebibliography}
\end{document}